\documentclass[a4paper,11pt]{ijms-preprint}

\usepackage{graphicx}

\newtheorem{example}{Example}[section]

\begin{document}


\title{A Discrete Algorithm to the Calculus of Variations\thanks{This 
is a preprint of a paper whose final and definite form appeared 
in \emph{Int. J. Math. Stat. 9 (2011), no.~S11, 26--41}.}}

\author{\textbf{C\'elia T. L. M. Pereira}$^\dag$\\
{\small \texttt{celia.pereira@estgoh.ipc.pt}}
\and \textbf{Pedro A. F. Cruz}$^\ddag$\\
{\small \texttt{pedrocruz@ua.pt}}
\and \textbf{Delfim F. M. Torres}$^\ddag$\\
{\small \texttt{delfim@ua.pt}}}

\date{$^\dag$College of Technology and Management
of Oliveira do Hospital\\
Polytechnic Institute of Coimbra\\
3400-124 Oliveira do Hospital, Portugal\\[0.3cm]
$^\ddag$Department of Mathematics\\
University of Aveiro\\
3810-193 Aveiro, Portugal}

\maketitle


\begin{abstract}
\noindent \emph{A numerical study of an algorithm proposed by Gusein Guseinov,
which determines approximations to the optimal solution of problems
of calculus of variations using two discretizations
and correspondent Euler-Lagrange equations, is investigated.
The results we obtain to discretizations of the brachistochrone problem
and Mani\`{a} example with Lavrentiev's phenomenon are compared
with the solutions found by other methods and solvers.
We conclude that Guseinov's method presents better solutions
in most of the cases studied.}

\medskip


\medskip


\end{abstract}


\section{Introduction}

In a problem of the calculus of variations,
functions that extremize a given functional are sought.
Several methods to determine approximated solutions for such problems
are known in the literature. Here we do a comparison study
of different approaches. More precisely, we investigate a method proposed
by Guseinov in \cite{Guseinov} that determines approximated solutions
to problems of the calculus of variations by discretization.
The results we obtain applying Guseinov's method to well known problems
are then compared to the results found by other methods.


\section{Calculus of variations}

Given real numbers $x_a$, $x_b$, $\alpha$, $\beta$,
$x_a < x_b$, and a function $f : [x_a,x_b] \times \mathbb{R} \times \mathbb{R} \rightarrow \mathbb{R}$
called the Lagrangian, $f \in C^1$,
the fundamental problem of the calculus of variations aims to determine
$y \in C^1([x_a,x_b])$ such that function $y$ minimizes the functional
\begin{equation}
\label{problema_cv}
\mathcal{J}[y(\cdot)]=\int_{x_a}^{x_b} {f\left( {x,y(x),y'(x)} \right)dx}
\end{equation}
under the given boundary conditions
\begin{equation}
\label{eq:bc}
y(x_a) = \alpha \, , \quad y(x_b) = \beta \, .
\end{equation}

A necessary condition for function $y$ to solve problem
\eqref{problema_cv}-\eqref{eq:bc} is given by the Euler-Lagrange equation
\begin{equation}
\label{eq:euler-lagrange}
\frac{\partial f}{\partial y}\left(x,y(x),y'(x)\right)
=\frac{d}{d x}\left(\frac{\partial f}{\partial y'}\left(x,y(x),y'(x)\right)\right) \, .
\end{equation}
A solution $y$ of \eqref{eq:euler-lagrange} is said to be an \emph{extremal},
meaning that $y$ either minimizes, maximizes,
or acts like a ``saddle function'' of $\mathcal{J}[\cdot]$.


\begin{example}[the brachistochrone problem]
Given two points $A=(x_a,y_a)$ and $B=(x_b,y_b)$ in the same vertical plane,
we want to determine the curve described by a particle that, with zero initial
speed and under action of gravity, connects the points in a minimum time,
ignoring friction:
\begin{equation}
\label{problema_braqui}
\displaystyle \min_{y(\cdot)} \frac{1}{\sqrt{2g}}\int_{x_a}^{x_b}
\sqrt{\frac{1+y'(x)^2}{y_a - y(x)}}dx \, ,
\end{equation}
where $g$ is the gravitational constant.

The solution to problem \eqref{problema_braqui} is an arc of cycloid:
the curve described by the revolution of a point belonging to a circle of radius $r$
that rolls without sliding on the straight line $y = y_a$.
Such curve is described by the parametrization
\begin{equation*}
\left\{
\begin{array}{l}
x=x_a+r\left(\theta - \sin(\theta)\right)\\
y=y_a-r\left(1-\cos(\theta)\right)
\end{array}
,\quad \theta_0 \leq \theta \leq \theta_1.
\right.
\end{equation*}
Since we know $(x_a,y_a)$ and $(x_b,y_b)$, it is possible to determine $r$,
$\theta_0$ and $\theta_1$. Using the point $(x_a,y_a)$ and considering
$r \neq 0$, $\theta_0=0$ is obtained.
The time it takes the particle to travel the curve is given by
\begin{equation}
\label{funcional_continuo_braqui}
\frac{1}{\sqrt{2g}}\int_{x_a}^{x_b}\sqrt{\frac{1+y'(x)^2}{y_a - y(x)}}dx
= \frac{1}{\sqrt{2g}}\int_{0}^{\theta_1}\sqrt{\frac{1+\left(\frac{\frac{d y}{d
\theta}(\theta)}{\frac{d x}{d \theta}(\theta)}\right)^2}{y_a - y(\theta)}}\frac{d x}{d \theta}(\theta) d \theta.
\end{equation}
For instance, the brachistochrone problem \eqref{problema_braqui} with boundary conditions
\begin{eqnarray}
\label{condicoes_fronteira_braqui_expl}
x_a = 0\, , \quad y_a = 10\, , \quad x_b = 10\, , \quad y_b = 0
\end{eqnarray}
has the following function as optimal solution:
$$\left\{
\begin{array}{l}
x = 5.729170\theta-5.729170\sin(\theta)\\
y = 4.270830+5.729170\cos(\theta)
\end{array}
,\quad 0 \leq \theta \leq 2.412011.\right.$$

\begin{figure}[!ht]
\centering
\includegraphics[width=8cm]{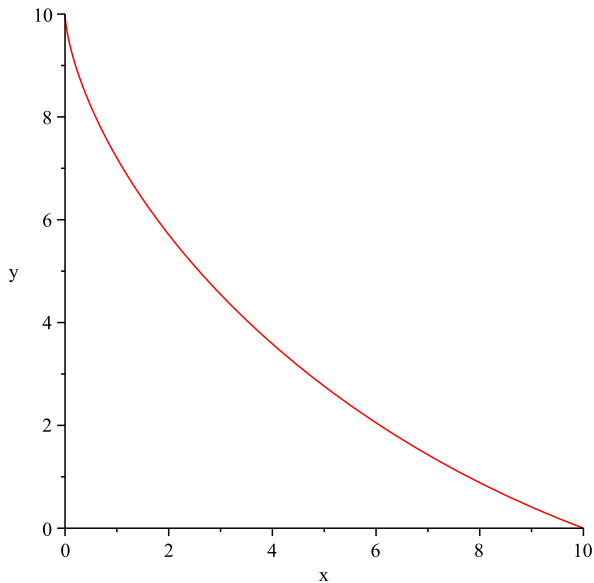}
\caption{Solution to the brachistochrone problem with $A=(0,10)$ and $B=(10,0)$.}
\end{figure}

Considering \eqref{funcional_continuo_braqui},
the optimal value, disregarding the constant $\frac{1}{\sqrt{2g}}$, is
\begin{eqnarray*}
\int_0^{2.412011}\sqrt{\frac{1+\left(
\frac{32.823393\sin(\theta)}{5.729170-5.729170\cos(\theta)}\right)^2}{5.729170
-5.729170\cos(\theta)}}\left(5.729170-5.729170\cos(\theta)\right)d\theta
= 8.164699.
\end{eqnarray*}
\end{example}


\begin{example}[Mani\`a's example]
The problem consists of finding the function $y$ that minimizes
\begin{equation}
\label{int_expl_mania}
\mathcal{I}[y(\cdot)] = \int^{1}_{0}(y(x)^3-x)^2(y'(x))^6dx
\end{equation}
under the boundary conditions $y(0)=0$ and $y(1)=1$.
Mani\`a example is a problem of the calculus of variations that exhibits
the so called Lavrentiev phenomenon \cite{Cesari,Ferriero}.
Problems with the Lavrentiev phenomenon have different solutions
in the space of absolutely continuous functions and in $C^1$ or the space of Lipschitzian functions:
the infimum of the functional in the set of absolutely continuous functions
is strictly less than the infimum of the same functional in the admissible set $C^1$.
This gap makes it harder (or even impossible) to determine the minimum
of the variational functional computationally.
The optimal solution for the Mani\`a problem is the function $\hat{y}(x)=x^\frac{1}{3}$.
Indeed, for all $x \in [0,1]$ one has
$L\left(x, y(x),y'(x)\right) = \left(y(x)^3-x\right)^2\left(y'(x)\right)^6 \geq 0$,
and thus $\mathcal{I}[y] \geq 0$. Since
$L\left(x, \hat{y}(x),\hat{y}'(x)\right)
= \left(\hat{y}(x)^3-x\right)^2\left(\hat{y}'(x)\right)^6 = 0$
for all $x \in [0,1]$, then $\mathcal{I}[\hat{y}(\cdot)] = 0$
is the optimum value for \eqref{int_expl_mania}.
Although the optimal solution for this problem is simple,
it is an open question how to find the optimal value zero for
$\mathcal{I}$ by means of numerical methods.
\end{example}


\section{Discrete-time calculus of variations}

A way to find approximate solutions to problems
\eqref{problema_cv}-\eqref{eq:bc} consists
in subdividing the problem into ``smaller'' problems easier to solve.
With this in mind, the interval $[x_a,x_b]$ is considered as a union of $n$ intervals:
\begin{equation*}
[x_a,x_b]=[x_a,x_1]\cup[x_1,x_2]\cup\cdots\cup[x_{n-1},x_b] \, .
\end{equation*}
Considering all the subintervals with the same amplitude,
\begin{equation}
\label{eq:interv}
[x_a,x_b]= \displaystyle \bigcup_{i=1}^n[x_{i-1},x_i] \, ,
\end{equation}
$x_0 = x_a$, $x_n = x_b$, $x_i=x_0+i\Delta x$, $i=1,2,\ldots,n-1$,
and $\Delta x = \frac{x_b-x_a}{n}$, $y'$ is approximated by the slope
of the line defined by the two extreme values of each interval:
$$y'(x_i) \approx \frac{y(x_i)-y(x_{i-1})}{\Delta x}
=\frac{y_i-y_{i-1}}{\Delta x}=\frac{\Delta y_{i-1}}{\Delta x}.$$


\subsection{The standard discretization}

There are several ways to formulate the problem using the intervals \eqref{eq:interv}.
The most common approach was presented by Euler himself, and we call it the ``standard''
discretization, also known as Euler's method of finite differences.
In this standard approach one searches $y_i$, $i=1,2,\ldots,n-1$, that minimize
(or maximize) the finite sum
\begin{equation*}
\displaystyle \sum^n_{i=1}f\left(x_i, y_i,\frac{\Delta y_{i-1}}{\Delta x}\right)\Delta x.
\end{equation*}
The Euler-Lagrange equations for this formulation are
\begin{equation*}
\frac{\partial f}{\partial y} \left(x_i, y_i, \frac{\Delta y_{i-1}}{\Delta x}\right) =
 \frac{\Delta \frac{\partial f}{\partial y'} \left(x_i,y_i,\frac{\Delta y_{i-1}}{\Delta x}\right)}
 {\Delta x},
\end{equation*}
$i=1,2,\ldots,n-1$.


\begin{example}[brachistochrone problem]
Using the standard discretization, the brachistochrone problem is approximated as follows:
$$
\begin{array}{cl}
\min & \displaystyle \sum_{i=1}^n \sqrt{\frac{1+\left(\frac{\Delta
y_{i-1}}{\Delta x_{i-1}}\right)^2}{y_a - y_i}}\Delta x_{i-1}\\
\text{subject to } & y(x_a) = y_a \, , \quad  y(x_b) = y_b.
\end{array}
$$
\end{example}


\begin{example}[Mani\`a's example]
The standard discretization of the Mani\`a example is given by
$$
\begin{array}{cl}
\min & \displaystyle \sum_{i=1}^n (y_i^3-x_i)^2\left(\frac{\Delta
y_{i-1}}{\Delta x_{i-1}}\right)^6\Delta x_{i-1}\\
\text{subject to } & y(0) = 0\, , \quad y(1) = 1.
\end{array}
$$
\end{example}


\subsection{Guseinov's discretization}

Guseinov proposes in \cite{Guseinov} a different discrete approach to find solutions to problems of the calculus of variations.
Considering a finite set of $n+1$ real numbers from the integration interval,
$$X = \{x_0, x_1, \ldots, x_n\} \, ,$$
$x_a=x_0 < x_1 < \cdots <x_n=x_b$, another set of $n+1$ real numbers
is found:
$$\{y_0=y(x_0)=\alpha, y_1=y(x_1), \ldots, y_n=y(x_n)=\beta\}.$$
Then, it is possible to define a polygonal line that connects all the points
$$(x_0,\alpha),(x_i,y_i),(x_n,\beta)\, , \quad i=1,\ldots,n-1 \, .$$

Guseinov's discretization is different in the way the functional is determined.
The discrete functional used by Guseinov \cite{Guseinov} is obtained considering
the meaning of the problem and its implications when we discretize it. For instance,
if we search a curve that extremizes a functional, such line is approximated
using several straight line segments. If an area is searched, then we use trapezoids
to define the discrete functional; when we search a revolution solid
we approximate it with frustums of cones.
The discrete Lagrangian of Guseinov is defined by
$$
\begin{array}{lccl}
L:& X \times X \times \mathbb{R} \times \mathbb{R} \times \mathbb{R}
& \longrightarrow & \mathbb{R}\\
&(s,t,u,v,w)& \mapsto & L(s,t,u,v,w),\\
\end{array}
$$
where
$\frac{\partial L}{\partial u}$, $\frac{\partial L}{\partial v}$,
$\frac{\partial L}{\partial w} \in C^0[\mathbb{R}]$
and the variables $s$ and $t$ concern the domain
of the functional \eqref{problema_cv}, $u = y(s)$, $v = y(t)$,
and $w = \frac{v-u}{t-s}$.
The finite sum to minimize (or maximize) is given by
\begin{equation}
\label{funcional_discreto_guseinov}
\mathcal{J}[y] = \displaystyle
\sum^n_{i=1}L\left(x_{i-1},x_i,y_{i-1},y_i,\frac{\Delta y_{i-1}}{\Delta x_{i-1}}\right)\Delta x_{i-1},
\end{equation}
where $\Delta y_{i-1}=y_i - y_{i-1}$ and $\Delta x_{i-1}=x_i - x_{i-1}$.
The function $y$ to be found verify the boundary conditions and should minimize
(or maximize) $\mathcal{J}$.


\begin{example}[brachistochrone problem]
To discretize the brachistochrone problem \emph{a la} Guseinov \cite{Guseinov},
we consider that the solution curve will be approximated by a sequence
of $n$ straight line segments that connect the points $\left(x_{i-1},y_{i-1}\right)$
and $\left(x_i,y_i\right)$, $i=1,\ldots, n$, $x_i > x_{i-1}$, $y_{i-1} > y_i$.
The amount of time spent by a particle to travel each of these line segments
(neglecting friction) should be calculated.
Regarding $a_i$ as the acceleration of the particle and $\theta_i$
as the angle formed with the $x$'s axis of the segment $i$,
$a_i = g \sin \theta_i$. By simple trigonometry,
$$\sin \theta_i = \frac{y_{i-1}-y_i}{\sqrt{\left(x_i - x_{i-1}\right)^2
+\left(y_{i-1} - y_i\right)^2}}=
-\frac{\Delta y_{i-1}}{\sqrt{\left(\Delta x_{i-1}\right)^2
+ \left(-\Delta y_{i-1}\right)^2}}.$$
Thus, $a_i = -\frac{g\Delta y_{i-1}}{\sqrt{\left(\Delta x_{i-1}\right)^2
+ \left(\Delta y_{i-1}\right)^2}}$.
Because $a_i = \frac{dv_i}{dt}$, $\frac{dv_i}{dt}
=-\frac{g\Delta y_{i-1}}{\sqrt{\left(\Delta x_{i-1}\right)^2
+ \left(\Delta y_{i-1}\right)^2}}$, then
$$v_i=-\frac{g\Delta y_{i-1}}{\sqrt{\left(\Delta x_{i-1}\right)^2
+ \left(\Delta y_{i-1}\right)^2}}t + c.$$
Besides, if $t=0$, then $c=v_0$. In each line segment
the counting of time is resumed (we count the time needed
to run each line segment separately, adding all the times at the end).
This way, in each line segment,
$t_{\mbox{initial}}=0$. Thus, $c=v_{i-1}$, $v_0 = 0$, and
$$
v_i=-\frac{g\Delta y_{i-1}}{\sqrt{\left(\Delta x_{i-1}\right)^2
+ \left(\Delta y_{i-1}\right)^2}}t+v_{i-1}.
$$
Moreover, $v = \frac{ds}{dt}$, and therefore
$$\displaystyle \int_0^{t_i} v dt = \int ds,$$
where $t_i$ is the time the particle takes to run the line segment $i$:
\begin{equation*}
\begin{split}
\displaystyle\int_0^{t_i} &\left(-\frac{g \Delta y_{i-1}}{\sqrt{\left(\Delta
x_{i-1}\right)^2 + \left(\Delta y_{i-1}\right)^2}}t + v_{i-1}\right)dt\\
&\qquad =-\frac{g\Delta y_{i-1}}{\sqrt{\left(\Delta x_{i-1}\right)^2
+ \left(\Delta y_{i-1}\right)^2}}\displaystyle\int_0^{t_i} t dt
+ v_{i-1}\displaystyle\int_0^{t_i} dt\\
&\qquad =-\frac{g\Delta y_{i-1}}{2\sqrt{\left(\Delta x_{i-1}\right)^2
+ \left(\Delta y_{i-1}\right)^2}} t_i^2 + v_{i-1}t_i\\
&\qquad =-\frac{g\Delta y_{i-1}}{2\Delta x_{i-1}\sqrt{1
+ \left(\frac{\Delta y_{i-1}}{\Delta x_{i-1}}\right)^2}} t_i^2 + v_{i-1}t_i \, .
\end{split}
\end{equation*}
To determine the line integral $ds$, a parametrization of the line segment
which connects the points $\left(x_{i-1},y_{i-1}\right)$ and $\left(x_i,y_i\right)$ is used:
$$
\begin{cases}
x(t) = t \, ,\\
y(t) = \frac{\Delta y_{i-1}}{\Delta x_{i-1}}t+y_{i-1}
-\frac{\Delta y_{i-1}}{\Delta x_{i-1}}x_{i-1} \, ,
\end{cases}
$$
$t \in [x_{i-1},x_i]$, and
$$
\begin{array}{lcl}
\displaystyle \int ds &=& \displaystyle \int_{x_{i-1}}^{x_i}
\sqrt{\left(\frac{d x(t)}{d t}\right)^2 + \left(\frac{d y(t)}{d t}\right)^2}dt\\
 &=& \displaystyle \int_{x_{i-1}}^{x_i} \sqrt{1
 + \left(\frac{\Delta y_{i-1}}{\Delta x_{i-1}}\right)^2}dt\\
 &=& \sqrt{1 + \left(\frac{\Delta y_{i-1}}{\Delta x_{i-1}}\right)^2}\int_{x_{i-1}}^{x_i}dt\\
 &=& \sqrt{1 + \left(\frac{\Delta y_{i-1}}{\Delta x_{i-1}}\right)^2}\left(x_i -x_{i-1}\right)\\
 &=& \Delta x_{i-1}\sqrt{1 + \left(\frac{\Delta y_{i-1}}{\Delta x_{i-1}}\right)^2}.
\end{array}
$$
Then,
\begin{multline*}
-\frac{g\Delta y_{i-1}}{2\Delta x_{i-1}\sqrt{1 +
\left(\frac{\Delta y_{i-1}}{\Delta x_{i-1}}\right)^2}} t_i^2 + v_{i-1}t_i
= \Delta x_{i-1}\sqrt{1 + \left(\frac{\Delta y_{i-1}}{\Delta x_{i-1}}\right)^2}\\
\Leftrightarrow  \frac{g\Delta y_{i-1}}{2\Delta x_{i-1}\sqrt{1
+ \left(\frac{\Delta y_{i-1}}{\Delta x_{i-1}}\right)^2}} t_i^2 - v_{i-1}t_i
+ \Delta x_{i-1}\sqrt{1 + \left(\frac{\Delta y_{i-1}}{\Delta x_{i-1}}\right)^2} = 0\, .
\end{multline*}
Using the quadratic formula,
\begin{equation}
\label{eq:ti}
t_i = \frac{\Delta x_{i-1}\sqrt{1+\left(\frac{\Delta y_{i-1}}{\Delta
x_{i-1}}\right)^2}}{g\Delta y_{i-1}}\left(v_{i-1} \pm \sqrt{v_{i-1}^2 -2g\Delta y_{i-1}}\right).
\end{equation}
Since $\Delta y_{i-1} < 0$ and $\Delta x_{i-1} > 0$,
then $\frac{\Delta x_{i-1}\sqrt{1+\left(\frac{\Delta y_{i-1}}{\Delta x_{i-1}}
\right)^2}}{g\Delta y_{i-1}} < 0$. Moreover, since the amount
of time must be non-negative, the condition
$$v_{i-1} \pm \sqrt{v_{i-1}^2 - 2g\Delta y_{i-1}} \leq 0$$
must be true and
$$
\begin{array}{rl}
& v_{i-1}^2 \geq v_{i-1}^2 - 2g\Delta y_{i-1}
\Leftrightarrow \left|v_{i-1}\right| \geq \sqrt{v_{i-1}^2 - 2g\Delta y_{i-1}}\\
\Rightarrow & v_{i-1} \geq \sqrt{v_{i-1}^2 - 2g\Delta y_{i-1}}\quad \text{ (since $v_{i-1} \geq 0$)}\\
\Rightarrow & v_{i-1} - \sqrt{v_{i-1}^2 - 2g\Delta y_{i-1}} \geq 0 \, .
\end{array}
$$
Since the expression should depend only on the points used to discretize the curve,
the initial velocity of the particle in each line segment $i$ ($v_{i-1}$)
shall be determined using these values. Applying the energy conservation law
(for the importance of conservation laws in the calculus of variations
we refer the reader to \cite{Gouveia1,Gouveia2}), the amount of mechanical energy in the beginning
of the first straight line segment ($E_{m_a}$) is the same as the amount
of energy at the beginning of each of the other line segments. Considering the line segment $i$,
$E_{m_a} = E_{m_i}$. At any given moment, the mechanical energy of the particle is the sum
of its potential and kinetic energies.
Considering $g$ as the gravity acceleration,
$m$ the mass of the particle, $h$ the height at which the particle is,
and $v$ the velocity of the particle, then
$$E_p = gmh \quad \text{and} \quad E_c = \frac{mv^2}{2}.$$
Therefore, in the line segment $i$ (with extremes
at $(x_{i-1},y_{i-1})$ and $(x_i,y_i)$),
\begin{equation*}
\begin{split}
E_{m_a} = E_{m_i} \Leftrightarrow & E_{p_a} + E_{c_a} = E_{p_i} + E_{m_i}\\
\Leftrightarrow & gmy_a + \frac{m \displaystyle v^2_0}{2}
= gmy_{i-1} + \frac{m \displaystyle v^2_{i-1}}{2}\\
\Leftrightarrow & gy_a + \frac{\displaystyle v^2_0}{2}
= gy_{i-1} + \frac{\displaystyle v^2_{i-1}}{2}.
\end{split}
\end{equation*}
Since the particle starts still, \textrm{i.e.}, $v_0 = 0$,
$$
gy_a = gy_{i-1} + \frac{\displaystyle v^2_{i-1}}{2}
\Leftrightarrow v_{i-1} = \pm \sqrt{2g\left(y_a-y_{i-1}\right)}.
$$
Besides, $v_{i-1} \geq 0$ for any $i=1, 2, \ldots, n$, so
$v_{i-1} = \sqrt{2g\left(y_a-y_{i-1}\right)}$.
Replacing this last expression for $v_{i-1}$ into \eqref{eq:ti},
the amount of time that a particle takes to run the line segment $i$ is
determined:
\begin{eqnarray*}
t_i & = & \frac{\Delta x_{i-1}\sqrt{1+\left(\frac{\Delta y_{i-1}}{\Delta
x_{i-1}}\right)^2}}{g\Delta y_{i-1}}\left(\sqrt{2g\left(y_a-y_{i-1}\right)}
- \sqrt{2g\left(y_a-y_{i-1}\right) - 2g\left(y_i-y_{i-1}\right)}\right)\\
& = & \frac{\Delta x_{i-1}\sqrt{1+\left(\frac{\Delta y_{i-1}}{\Delta
x_{i-1}}\right)^2}}{g\Delta y_{i-1}}\left(\sqrt{2g\left(y_a-y_{i-1}\right)}
- \sqrt{2g\left(y_a-y_i\right)}\right)\\
& = & \frac{\sqrt{2g}\sqrt{1+\left(\frac{\Delta y_{i-1}}{\Delta
x_{i-1}}\right)^2}}{g\Delta y_{i-1}}\left(\sqrt{y_a-y_{i-1}}
- \sqrt{y_a-y_i}\right)\Delta x_{i-1}\\
& = & \frac{\sqrt{2}\sqrt{1+\left(\frac{\Delta y_{i-1}}{\Delta
x_{i-1}}\right)^2}\left(\sqrt{y_a-y_{i-1}} - \sqrt{y_a-y_i}\right)
\cdot\left(\sqrt{y_a-y_{i-1}} + \sqrt{y_a-y_i}\right)\Delta x_{i-1}}{\sqrt{g}\Delta y_{i-1}\left(\sqrt{y_a-y_{i-1}}
+ \sqrt{y_a-y_i}\right)}\\
& = & \sqrt{\frac{2}{g}} \frac{\sqrt{1+\left(\frac{\Delta y_{i-1}}{\Delta
x_{i-1}}\right)^2}}{\Delta y_{i-1}\left(\sqrt{y_a-y_{i-1}}
+ \sqrt{y_a-y_i}\right)}\left(y_a - y_{i-1} - y_a + y_i\right)\Delta x_{i-1}\\
\end{eqnarray*}
\begin{eqnarray*}
& = & \sqrt{\frac{2}{g}} \frac{\sqrt{1+\left(\frac{\Delta y_{i-1}}{\Delta
x_{i-1}}\right)^2}}{\Delta y_{i-1}\left(\sqrt{y_a-y_{i-1}} +
\sqrt{y_a-y_i}\right)}\left(- y_{i-1} + y_i \right)\Delta x_{i-1}\\
& = & \sqrt{\frac{2}{g}} \frac{\sqrt{1+\left(\frac{\Delta y_{i-1}}{\Delta
x_{i-1}}\right)^2}}{\sqrt{y_a-y_{i-1}} + \sqrt{y_a-y_i}}\Delta x_{i-1}.
\end{eqnarray*}
Adding the times of the $n$ line segments we obtain
$$
\displaystyle \sum_{i=1}^n t_i = \sqrt{\frac{2}{g}}\displaystyle
\sum_{i=1}^n \frac{\sqrt{1+\left(\frac{\Delta y_{i-1}}{\Delta
x_{i-1}}\right)^2}}{\sqrt{y_a-y_{i-1}} + \sqrt{y_a-y_i}}\Delta x_{i-1}.
$$
To find an approximated solution for the brachistochrone problem, the $y_i$'s are searched so that
$$
\min \displaystyle \sum_{i=1}^n \frac{\sqrt{1+\left(\frac{\Delta
y_{i-1}}{\Delta x_{i-1}}\right)^2}}{\sqrt{y_a-y_{i-1}} + \sqrt{y_a-y_i}}\Delta x_{i-1}.
$$
The problem may, then, be formulated in the discrete (Guseinov) form as
$$
\begin{array}{cl}
\min & \displaystyle \sum_{i=1}^n \frac{\sqrt{1+\left(\frac{\Delta
y_{i-1}}{\Delta x_{i-1}}\right)^2}}{\sqrt{y_a-y_{i-1}} + \sqrt{y_a-y_i}}\Delta x_{i-1}\\
\text{subject to } & y(x_a) = y_a\, , \quad  y(x_b) = y_b.
\end{array}
$$
\end{example}


\begin{remark}
Since the Mani\`a example is a theoretical problem,
without no \emph{a priori} physical meaning, it is not clear
how to discretize it using Guseinov's approach.
\end{remark}


\subsection{Euler-Lagrange equation}

The next result presents the Euler-Lagrange equation
for the discrete Guseinov formulation of the problem
of the calculus of variations.

\begin{theorem}[Guseinov's discrete Euler-Lagrange equation \cite{Guseinov}]
\label{teorema-e-l-guseinov}
If the variational functional of the calculus of variations \eqref{funcional_discreto_guseinov}
has a local extreme in $\hat{y}$, then $\hat{y}$ satisfies the Euler-Lagrange equation
\begin{multline}
\label{eq:e-l-guseinov}
\frac{\partial L}{\partial u}\left(x_i,x_{i+1},y_i,y_{i+1},
\frac{\Delta y_i}{\Delta x_i}\right)\frac{\Delta x_i}{\Delta x_{i-1}}
+ \frac{\partial L}{\partial v}\left(x_{i-1},x_i,y_{i-1},y_i,\frac{\Delta y_{i-1}}{\Delta x_{i-1}}\right)\\
= \frac{\Delta \frac{\partial L}{\partial w}\left(x_{i-1},x_i,y_{i-1},y_i,
\frac{\Delta y_{i-1}}{\Delta x_{i-1}}\right)}{\Delta x_{i-1}}
\end{multline}
for $i \in \{1,2,\ldots,n-1\}$.
\end{theorem}

\begin{proof}
See \cite{Guseinov}.
\end{proof}

\begin{remark}
The Euler-Lagrange equation \eqref{eq:e-l-guseinov}
is to be complemented with the given
boundary conditions $\hat{y}(x_0)=\alpha$ and $\hat{y}(x_n)=\beta$.
\end{remark}

Using Theorem~\ref{teorema-e-l-guseinov} (and the given boundary conditions),
Guseinov \cite{Guseinov} developed a (new) method that, discretizing the interval,
determines approximated solutions for problems of the calculus of variations.
We note that Theorem~\ref{teorema-e-l-guseinov} does not imply the use of a
``Guseinov discretization''. There is no reference to the \emph{kind} of discrete functional $L$.
Functional $L$ must, only, verify the conditions $\frac{\partial L}{\partial u}$,
$\frac{\partial L}{\partial v}$, $\frac{\partial L}{\partial w} \in C^0[\mathbb{R}]$.
In other words, the method may be applied using the ``standard discretization''
as long as the standard discrete functional
verifies such conditions. The approximations to solutions
found by application of Theorem~\ref{teorema-e-l-guseinov}
are candidates to extremize the functional (they are \emph{critical functions}),
which means that the found functions may maximize, minimize or neither maximize nor minimize
the functional (\emph{saddle function}).


\subsection{Guseinov's algorithm}
\label{algimpl}

Based on Theorem~\ref{teorema-e-l-guseinov}, the next algorithm determines
approximated solutions for problems of the calculus of variations.

\bigskip

\textbf{Input:} a continuous Lagrangian $f$, a discrete Lagrangian $L$,
boundary conditions $y(x_a) = \alpha$ and $y(x_b) = \beta$,
and the number $n$ of intervals that divide the integration domain.\\

\textbf{Output:} $y(x_i)$, $i = 0,1,\ldots,n$, $\mathcal{J}[y]$
given by \eqref{funcional_discreto_guseinov},
and $\int_{x_a}^{x_b}f\left(x,\bar{y}(x),\bar{y}'(x)\right)dx$,
where $\bar{y}$ is the piecewise linear function defined by the points $y(x_i)$ found
by the algorithm.\\

\textbf{Algorithm}
\begin{enumerate}
\item Read input data: $L(s,t,u,v,w)$, $t_0 \leftarrow x_a$,
$t_n \leftarrow x_b$, $\alpha \leftarrow y_a$, $\beta \leftarrow y_b$, and $n$.
\item Determine the partial derivatives of $L$:
$Lu \leftarrow \frac{\partial L}{\partial u}\left(s,t,u,v,w\right)$,
$Lv \leftarrow \frac{\partial L}{\partial v}\left(s,t,u,v,w\right)$,
$Lw \leftarrow \frac{\partial L}{\partial w}\left(s,t,u,v,w\right)$.
\item Determine the set of $n+1$ equidistant $x$:
\begin{itemize}
\item $\Delta x=\frac{x_b-x_a}{n}$;
\item $X=\{x_0=x_a,x_1=x_0+\Delta x,\ldots,x_i=x_0+i\Delta x,\dots,x_n=x_b\}.$
\end{itemize}
\item Solve the discrete (nonlinear) Euler-Lagrange equations system \eqref{eq:e-l-guseinov}:
\begin{itemize}
\item determine the system of equations $F(Y) = 0$ such that
\begin{itemize}
\item for $i=1,\dots,n-1$,
$$
\begin{array}{ll}
F_i = & Lu\left(x_i,x_{i+1},y_i,y_{i+1},\frac{y_{i+1} -
y_i}{\Delta x}\right)
+ Lv \left(x_{i-1},x_i,y_{i-1},y_i,\frac{y_i - y_{i-1}}{\Delta x}\right)\\
& - \frac{Lw\left(x_i,x_{i+1},y_i,y_{i+1},\frac{y_{i+1} -
y_i}{\Delta x}\right) -
Lw\left(x_{i-1},x_i,y_{i-1},y_i,\frac{y_i - y_{i-1}}{\Delta x}\right)} {\Delta x}
\end{array}
$$
\item the boundary conditions
$F_0 = y_0 - \alpha$, $F_n = y_n - \beta$ hold;
\end{itemize}
\item solve the system.
\end{itemize}
\item Compute $\mathcal{J}[y(\cdot)] = \displaystyle \sum_{i=1}^{n}
L\left(x_{i-1},x_i,y(x_{i-1}),y(x_i),\frac{\Delta y_{i-1}}{\Delta x}\right)\Delta x$.
\item Determine the piecewise linear function $\bar{y}$.
\item Compute the value of $\int_{x_a}^{x_b}f(x,\bar{y}(x),\bar{y}'(x))dx$.
\item Draw the graphic of function $\bar{y}$.
\end{enumerate}


\subsection{Implementation}

The algorithm of \S\ref{algimpl} was implemented using \emph{Mathematica}~6.
Considering the problem and the possible amount of points chosen from the domain of the integral,
the nonlinear system of equations is not, in general, easily solved. So, in general, the \emph{Mathematica}'s
functions used to solve the system of equations may not find solutions or may find complex numbers.
Neither of these results is acceptable.
In our implementation, a function that numerically solves nonlinear programming problems is used.
In this way an approximation to the solution of the system of equations is found.
Given the nonlinear equations system
$$
\left\{
\begin{array}{l}
F(i)=0,\mbox{ }i=1,\ldots,n-1\\
y_0-\alpha=0\\
y_n-\beta=0
\end{array}
\right.
$$
we consider the following
nonlinear programming problem:
\begin{equation}
\label{non:prb}
\begin{array}{rl}
\min&\displaystyle \sum_{i=1}^{n-1}\left(F(i)\right)^2
+\left(y_0-\alpha\right)^2+\left(y_n-\beta\right)^2\\
\text{subject to } &y_0=\alpha\, , \quad y_n=\beta\, ,\\
&y_i \in \mathbb{R} \, .
\end{array}
\end{equation}
To find a numerical approximation to the solution of the nonlinear programming
problem \eqref{non:prb}, the function \emph{NMinimize} of \emph{Mathematica} is used.
It is possible to choose the method used to find solutions:
\emph{Random Search}, \emph{Nelder Mead}, \emph{Differential Evolution}, or \emph{Simulated Annealing}.
These methods present solutions that may or may not be the same. Besides, one of the solutions
may be the best considering the nonlinear problem (solving the nonlinear equations system)
but not the best regarding the initial problem. So, from the presented solutions by these methods,
the one that minimizes the discrete problem of the calculus of variations is the one chosen.
Although some methods end with success, they may present solutions that make no sense
(such as $y_i$ close to plus or minus infinity).
Thus, since the nonlinear programming problem may have several restrictions,
a new set of restrictions was added, trying to improve the solutions.
This way, two set of restrictions were tested:
\begin{itemize}
\item ``restrictions 1'': boundary conditions from the problem of the calculus of variations;
\item ``restrictions 2'': boundary conditions from the problem of the calculus of variations and
\begin{equation}
\label{restricoes2_mathematica}
|y_i| \leq \max\{|\alpha|,|\beta|\},\quad i=1,\ldots,n-1.
\end{equation}
\end{itemize}


\subsection{An optimal control solver}

\emph{OC} -- Optimal Control solver -- is a solver supplied in \cite{Smirnov}.
This solver determines approximated solutions for optimal control problems
in Mayer formulation. The solutions are found by discretizing the interval $[x_a,x_b]$
in $n$ subintervals and using Mathematical Programming methods.
It is possible to choose the method used to optimize
(\emph{Conjugated Gradients} -- CG, \emph{Newton Method} -- NM, \emph{Univariate Search} -- US,
\emph{Direct Search} -- DS 1 and DS 2, and \emph{Random Search} -- RS),
the method used to solve the differential equations
(\emph{Euler} or \emph{Runge-Kutta}), the initial solution,
maximum number of iterations, and the number of intervals of the discretization.
To use this solver, the problem of the calculus of variations must be formulated as an optimal control problem:
$$
\begin{array}{rl}
\min&\displaystyle\int_{x_a}^{x_b}{L\left( {x,y(x),z(x)} \right)dx}\\
\mbox{subject to } &y'(x)=z(x)\\
&y(x_a)=\alpha\, , \quad y(x_b)=\beta.
\end{array}
$$
This optimal control problem in Mayer form is:
$$
\begin{array}{cl}
\min & u(x_b)\\
\text{subject to } & y'(x)=z(x)\\
& u'(x)=L\left( {x,y(x),z(x)} \right)\\
& u(x_a)=0\\
& y(x_a) = \alpha\, , \quad y(x_b) = \beta.
\end{array}
$$


\subsection{An evolutionary algorithm}

A simplification of the $(\mu,\lambda)-$ES Algorithm (an evolutionary algorithm)
presented in \cite{Cruz_e_Torres}, uses evolutionary strategies combined with optimal control
to find approximated solutions to problems of the calculus of variations.
The algorithm keeps seeking solutions to the problem, evaluating each of them.
The solutions closer to the target set are used to find new solutions (supposedly better).
This process ends after a certain number of iterations (see \cite{Cruz_e_Torres} for details).


\section{Results and comparisons}

The result to be compared is the value of the
functional integral \eqref{problema_cv}
along the approximated solutions.


\begin{example}[brachistochrone problem]
Table~\ref{tab:tabela_braqui_teste1} presents the value of the integral
\eqref{funcional_continuo_braqui} using as integrand
the piecewise linear functions defined with the approximated solutions
of \eqref{problema_braqui} found by the Guseinov algorithm of \S\ref{algimpl}
and the piecewise linear function defined with the optimal solution (PLFOpt),
whose value is used as reference.

\begin{table}[!ht]
\caption{Values of $\int_{x_a}^{x_b}f(x,\bar{y}(x),\bar{y}'(x))dx$ for problem \eqref{problema_braqui}
subject to $x_a=0$, $x_b=10$, $y_a = 10$, and $y_b=0$.}
\centering
\begin{tabular}{|c||c|c|c|c|c|c|}
\hline
No. && \multicolumn{2}{|c|}{Guseinov discretization} & \multicolumn{2}{|c|}{Standard discretization} & Integral\\
\cline{3-7}
of & Method & Restrictions 1 & Restrictions 2 & Restrictions 1 & Restrictions 2 & PLFOpt\\
Ints && \eqref{condicoes_fronteira_braqui_expl} & \eqref{condicoes_fronteira_braqui_expl} and \eqref{restricoes2_mathematica}&&&\\
\hline
\hline
$3$& RS &8.353 & 8.353 & 8.418 & 8.418 & 8.369\\
& NM & 8.353 & 8.353 & 8.418 & 8.418&\\
& DE& * & * & * & 8.418&\\
& SA& 8.353 & 8.353 & 8.418 & 8.418&\\
\hline
\hline
$5$& RS & 8.271 & 8.271 & 8.464 & 8.464 & 8.281\\
& NM & 8.271 & 8.271 & 8.464 & 8.464&\\
& DE& * & * & * & *&\\
& SA& 8.271 & ** & 8.464 & 8.464&\\
\hline
\hline
$8$& RS & 8.229 & 8.229 & 8.566 & 8.566 & 8.235\\
& NM & 8.229 & 9.037 & 8.566 & 9.906&\\
& DE& * & * & * & *&\\
& SA & *** & 8.229 & 8.566 & 12.288&\\
\hline
\hline
$10$& RS & *** & 12.716 & 8.617 & 12.682 & 8.220\\
& NM & 8.216 & 12.716 & 8.617 & 8.617&\\
& DE & * & * & * & *&\\
& SA & *** & 12.716 & 8.617 & 15.370&\\
\hline
\hline
$15$& RS & *** & ** & 8.702 & 12.749 & 8.201\\
& NM & *** & 12.764 & 8.702 & 8.702&\\
& DE & * & * & * & *&\\
& SA & * & * & 8.702 & *&\\
\hline
\hline
$20$& RS & 8.189 & ** & 8.751 & 12.782 & 8.191\\
& NM & *** & 15.972 & 8.751 & 21.053&\\
& DE & * & 25.115 & * & *&\\
& SA & *** & * & * & *&\\
\hline
\end{tabular}
\begin{flushleft}
`*': method used to solve the nonlinear programming problem
doesn't end successfully (\textrm{e.g.}, the method exceeded the number of iterations);

`**': the value found is complex;

`***': the value found is too high and the solution makes no sense;

Methods: RS---\emph{Random Search}, NM---\emph{Nelder Mead}, DE---\emph{Differential Evolution},
SA---\emph{Simulated Annealing}.
\end{flushleft}
\label{tab:tabela_braqui_teste1}
\end{table}

The best result found by our implementation of the algorithm used
Guseinov's discretization, 20 intervals for discretization,
and the \emph{Random Search} method (Figure~\ref{fig:graf_braqui_gus1_20_rs}).
The integral value is $8.189$.
Using the same options  except the type of discretization -- standard discretization --
the value of the integral is $8.751$, which shows that the approximation is clearly worse.

\begin{figure}[!ht]
\centering
\includegraphics[width=8cm]{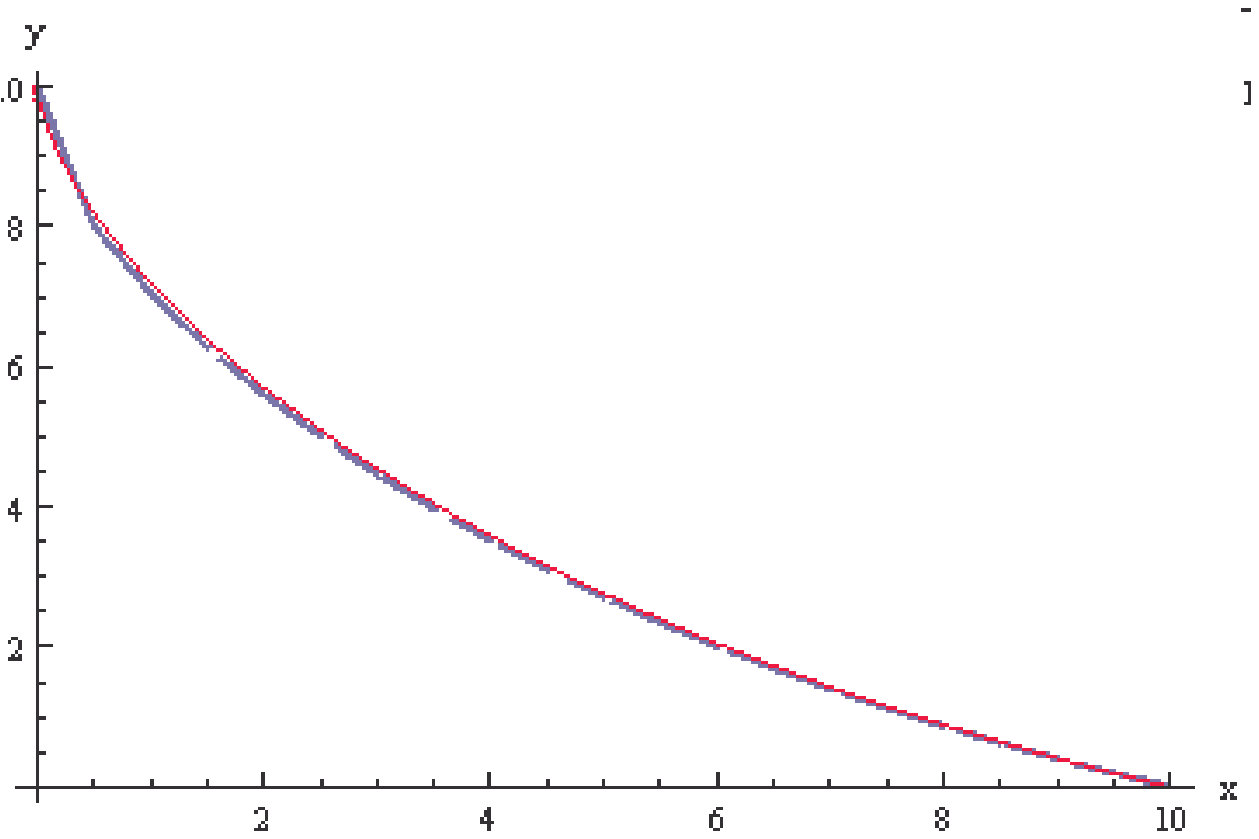}
\caption{Approximation found applying Guseinov's discretization
with $n=20$ intervals and the \emph{Random Search} method; and the exact solution.}
\label{fig:graf_braqui_gus1_20_rs}
\end{figure}

Considering the mentioned algorithms, solver,
and the solutions determined by them,
it is possible to verify which one presents
the best solution regarding the value of the integral
-- see Table~\ref{tab:sintese_braqui}.

\begin{table}[!ht]
\caption{\small{Values of the integral along approximated solutions obtained by different methods.}}
\centering
\begin{tabular}{|c|c|c|}
\hline
\small{Solution}& \small{Value of Integral} & \small{Notes} \\
\hline
\small{Optimal} & 8.16470 &\\
\hline
\small{Points from the optimal line} & 8.19139 & \small{21 points with}\\
&& \small{ (successive) equidistant abscissae} \\
\hline
\small{Guseinov algorithm} & 8.189344 & \small{20 intervals}\\
&&\small{Guseinov discretization}\\
&&\small{Method: \emph{Random Search}}\\
\hline
\small{Evolutionary algorithm} & 8.19365 & \small{20 intervals}\\
\hline
\small{\emph{OC}} & 8.336 & \small{20 intervals}\\
&& \small{Method: \emph{Conjugated Gradients}}\\
&& \small{Method: \emph{Runge-Kutta} 2}\\
&& \small{\emph{Piecewise Linear}}\\
\hline
\end{tabular}
\label{tab:sintese_braqui}
\end{table}
\end{example}

\bigskip


\begin{example}[Mani\`a's example]
The results obtained with our implementation
of Guseinov's algorithm (\S\ref{algimpl})
using the standard discretization with two different sets of restrictions
in the nonlinear programming problem and different methods are presented
on Table~\ref{tab:tabela_lavrentiev}.

\begin{table}[!ht]
\caption{Values of $\int_{x_a}^{x_b}L(x,\bar{y}(x),\bar{y}'(x))dx$
for the functional \eqref{int_expl_mania} ($x_a=0$, $x_b=1$, $y_a = 0$ and $y_b=1$).}
\centering
\begin{tabular}{|c||c|r|r|c|}
\hline
Number of & Method &\multicolumn{2}{|c|}{Standard discretization}&Integral\\
\cline{3-5}
 Intervals&& Restrictions 1 & Restrictions 2 & PLFOpt\\
\hline
\hline
$3$& RS &92.314&92.314&0.229\\
& NM & 92.314&92.314&\\
& DE&92.314&92.314&\\
& SA&92.314&92.314&\\
\hline
\hline
$5$& RS &1488.100&1488.100&0.381\\
& NM &1488.090&0.994&\\
& DE&1488.100&1488.100&\\
& SA &1488.100&1488.100&\\
\hline
\hline
$8$& RS &17549.400&14.768&0.610\\
& NM &539.515&537.198&\\
& DE&64.207&17549.400&\\
& SA&17549.400&0.328&\\
\hline
\hline
$10$& RS &55619.100&0.490&0.762\\
& NM &55619.100&20.364&\\
& DE &2390.020&55619.000&\\
& SA &55619.100&3.750&\\
\hline
\hline
$15$& RS &443732.000&1.156&1.143\\
& NM &443732.000&12416.200&\\
& DE &0.474&443732.000&\\
& SA &443732.000&66.265&\\
\hline
\hline
$20$& RS &1915810.000&127.346&1.524\\
& NM &1915810.000&10.174&\\
& DE &2235.890&32784.700&\\
& SA &0.956&7.473&\\
\hline
\end{tabular}
\label{tab:tabela_lavrentiev}
\end{table}

\bigskip

The solution of the best ``numerical'' result obtained using our implementation
in \emph{Mathematica} applied 8 discretization intervals and the method
of Simulated Annealing (Figure~\ref{fig:grafico_lavrentiev_alg_melhor_valor}).
The result is the value $0.328$ for the integral \eqref{int_expl_mania}.

\begin{figure}[!ht]
\centering
\includegraphics[width=7cm]{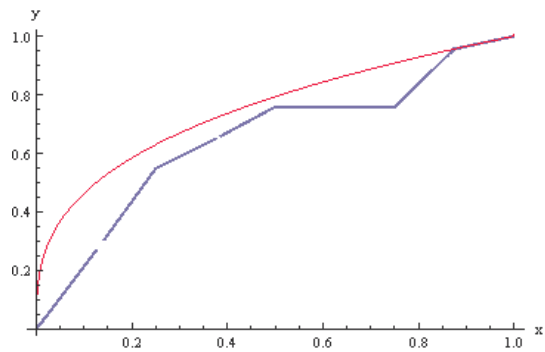}
\caption{\scriptsize{Approximation obtained by applying the standard discretization with
``Restrictions 2'', $n=8$, and the \emph{Simulated Annealing} method;
together with the optimal absolute continuous curve.}}
\label{fig:grafico_lavrentiev_alg_melhor_valor}
\end{figure}

Approximations for the solution found with $n=10$, $n=15$ or $n=20$ discretization intervals
seem to be closer to the optimal curve than de previous one
(Figure~\ref{fig:grafico_lavrentiev_alg_mais_prox_curva}).
However, the results with $n=10$, $n=15$ or $n=20$
are worse than with $n=8$ (for instance, the value of the functional
obtained with 15 intervals is $1.156$).

\begin{figure}[!ht]
\centering
\includegraphics[width=7cm]{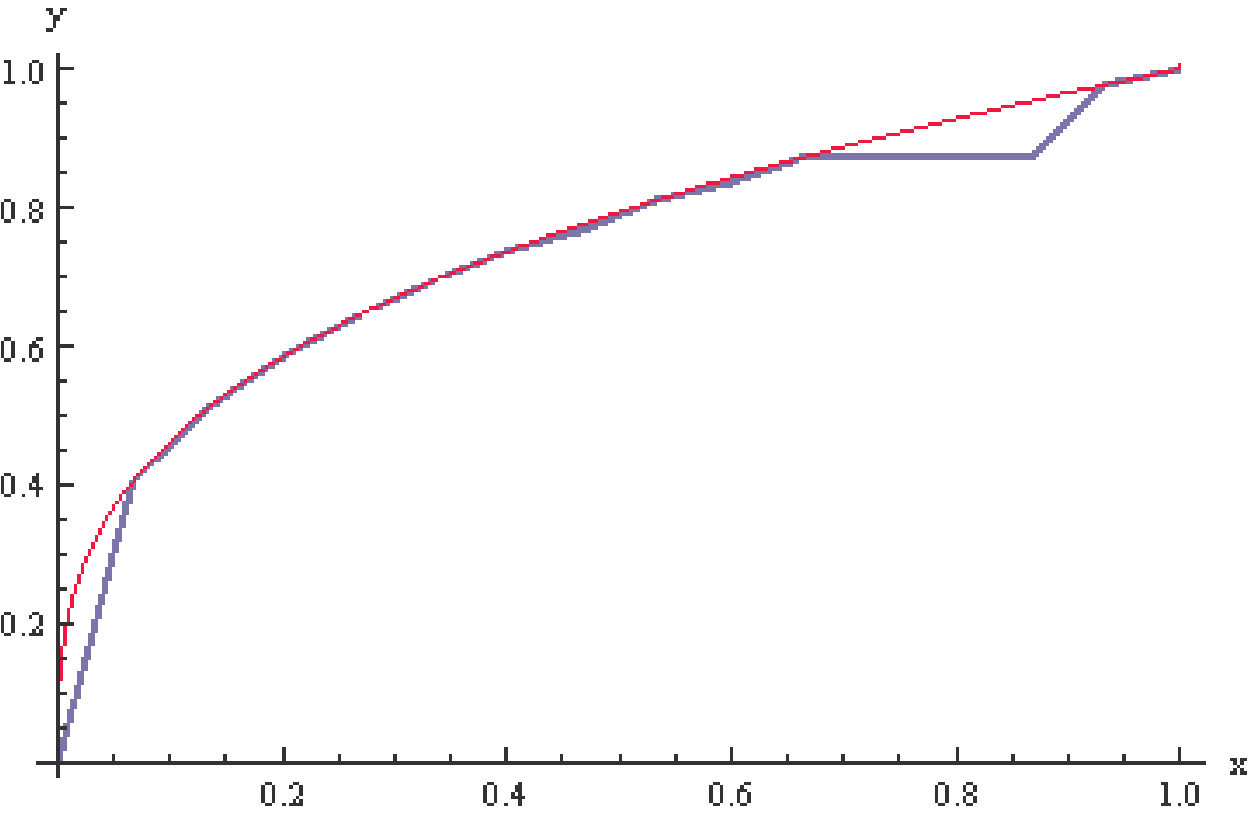}
\caption{\scriptsize{Approximation obtained by applying the standard
discretization with ``Restrictions 2'', $n=15$,
and the \emph{Random Search} method.}}
\label{fig:grafico_lavrentiev_alg_mais_prox_curva}
\end{figure}

Moreover, taking a close look at the last column of the Table~\ref{tab:tabela_lavrentiev},
something is apparently wrong because the value of the integral is increasing as the number
of discretization intervals increases, instead of becoming closer to zero (the optimal value).
Although this looks contradictory, it is not.  In fact this is a consequence of the Lavrentiev phenomenon
exhibited by this problem. The next result proposed and proved in \cite{Ferriero} explains the fact.

\begin{theorem}
\label{teor_f_lipschitzianas}
For any sequence of Lipschitz trajectories $\left\{y_n\right\}_n$
such that $y_n$ tends to $\hat{y}(x)=x^\frac{1}{3}$ as $n$ tends to $\infty$,
for almost all $x \in [0,1]$, then $\mathcal{I}[y_n]$
defined in \eqref{int_expl_mania} tends to $\infty$.
\end{theorem}

Considering the peculiarity of the Mani\`a example,
it is very difficult to compare objectively
the results determined by other methods and algorithms.
Unlike the brachistochrone problem, there is apparently no relation between
the solution that is graphically better and the solution that is numerically better.
However, regarding the graphics, the solutions found by the Guseinov algorithm are good.
Increasing the number of intervals used in the discretization, the approximations improve
considering the graphical representation but the integral value becomes worse.
This fact is also seen in the graphics of the piecewise linear functions defined using points
of the optimal solution and in the approximations found by the \emph{OC} solver.
The best approximation found by the \emph{OC} solver, considering the value of the integral,
used the methods of \emph{Conjugated Gradients}, \emph{Runge-Kutta}, and \emph{Piecewise Constant},
with 0.0326998 as the value of the integral. However the graphic of the approximation found using the methods
\emph{Univariate Search}, \emph{Runge-Kutta}, and \emph{Piecewise Linear},
shows that this approximation is closer to the optimal solution,
although its integral value is 1.45193.

In \cite{Bai_e_Li} the Truncation Method is proposed. This method determines an upper
limit for the integral minimum using an auxiliary functional whose integration domain
is in the original integration domain and does not include the points that
``create'' the Lavrentiev phenomenon. This new integration domain is also divided in $n$ intervals.
The numerical results presented in the paper \cite{Bai_e_Li}
don't include the integral value used in our work as the key comparison element.
\end{example}


\section{Conclusion}

The algorithm proposed by Guseinov in \cite{Guseinov}
presents very good solutions to ``regular'' problems of the
calculus of variations. Moreover, the usage of Guseinov discretization
improves the results. The method works worse when applied
to problems with the Lavrentiev gap, like the Mani\`a example, or
to problems of optimal control.
Since the method involves the resolution of nonlinear systems of equations,
if the numerical methods used in the solvers are not working well
the algorithm does not find solutions as good as it could.
The solutions obtained for the Mani\`a example are not worse
than the solutions found by other methods and solvers.
Our implementation of the algorithm in \emph{Mathematica} (version~6)
is very easy to use.


\section*{Acknowledgment}

Work partially supported by the R\&D unit
``Centre for Research on Optimization and Control''
of the University of Aveiro,
cofinanced by the European Community Fund FEDER/POCI 2010.



\end{document}